\documentclass[final,5p,times,twocolumn,authoryear]{elsarticle}

\usepackage{amsmath}
\usepackage{amssymb}
\usepackage{amsthm}
\usepackage{siunitx}

\DeclareMathOperator*{\maximize}{maximize}
\def\S{\mathbb{S}}
\def\X{\mathbb{X}}
\def\U{\mathbb{U}}
\def\W{\mathbb{W}}

\allowdisplaybreaks[4]

\begin{document}

\begin{frontmatter}

%% Title, authors and addresses

%% use the tnoteref command within \title for footnotes;
%% use the tnotetext command for theassociated footnote;
%% use the fnref command within \author or \affiliation for footnotes;
%% use the fntext command for theassociated footnote;
%% use the corref command within \author for corresponding author footnotes;
%% use the cortext command for theassociated footnote;
%% use the ead command for the email address,
%% and the form \ead[url] for the home page:
%% \title{Title\tnoteref{label1}}
%% \tnotetext[label1]{}
%% \author{Name\corref{cor1}\fnref{label2}}
%% \ead{email address}
%% \ead[url]{home page}
%% \fntext[label2]{}
%% \cortext[cor1]{}
%% \affiliation{organization={},
%%            addressline={}, 
%%            city={},
%%            postcode={}, 
%%            state={},
%%            country={}}
%% \fntext[label3]{}

\title{Deep learning adaptive Model Predictive Control of Fed-Batch Cultivations}
% \title{Data-driven Model Predictive Control of Fed-Batch Cultivations}

%% use optional labels to link authors explicitly to addresses:
%% \author[label1,label2]{}
%% \affiliation[label1]{organization={},
%%             addressline={},
%%             city={},
%%             postcode={},
%%             state={},
%%             country={}}
%%
%% \affiliation[label2]{organization={},
%%             addressline={},
%%             city={},
%%             postcode={},
%%             state={},
%%             country={}}

\author[2]{Niels Krausch}
\ead{niels.krausch@novartis.com}
\author[1]{Martin Doff-Sotta}
\ead{martin.doff-sotta@eng.ox.ac.uk}
\author[1]{Mark Cannon\corref{cor1}}
\ead{mark.cannon@eng.ox.ac.uk}
\author[2]{Peter Neubauer}
\ead{peter.neubauer@tu-berlin.de}
\author[2]{Mariano Nicolas Cruz Bournazou}
\ead{mariano.n.cruzbournazou@tu-berlin.de}

% \affiliation[1]{organization={Novartis Pharma AG},%Department and Organization
%             addressline={Werk Klybeck Klybeckstrasse 141}, 
%             city={Basel},
%             postcode={CH-4057}, 
%             country={Switzerland}}

\affiliation[1]{organization={University of Oxford, Department of Engineering Science},%Department and Organization
            addressline={Parks Road}, 
            city={Oxford},
            postcode={OX1 3PJ}, 
            country={UK}}

\affiliation[2]{organization={Technische Universität Berlin, Institute of Biotechnology, Bioprocess engineering},%Department and Organization
            addressline={Ackerstr.~76}, 
            city={Berlin},
            postcode={13355}, 
            country={Germany}}

\cortext[cor1]{Corresponding author}
          
\begin{abstract}
%% Text of abstract
  % Bioprocesses are often characterised by nonlinear and uncertain dynamics. This poses particular challenges in the context of model predictive control (MPC) algorithms because of their computational demand when applied to nonlinear systems.
  % %
  % Recent advances in optimal control theory have shown that concepts from convex optimisation, tube MPC, and differences of convex functions (DC) enable efficient robust online process control.
  % %
  % Our approach is based on DC decompositions of nonlinear dynamics and successive linearisations around predicted trajectories. By convexity, the linearisation errors have tight bounds and can be treated as bounded disturbances in a robust tube MPC framework.
  % %
  % We describe a systematic data-driven method of computing DC model representations using deep learning neural networks with a special convex structure and we explain how the resulting MPC optimisation can be solved using convex programming. 
  % % 
  % For the problem of maximising product formation in a cultivation with uncertain model parameters, we design a controller that ensures robust constraint satisfaction and allows unknown parameters to be estimated online.
  % %
  % Our results show that this is a promising approach for computationally tractable robust MPC of bioprocesses.
  %
  Bioprocesses are often characterised by nonlinear and uncertain dynamics, posing particular challenges for model predictive control (MPC) algorithms due to their computational demands when applied to nonlinear systems.
  Recent advances in optimal control theory have demonstrated that concepts from convex optimisation, tube MPC, and differences of convex functions (DC) enable efficient, robust online process control.
  Our approach is based on DC decompositions of nonlinear dynamics and successive linearisations around predicted trajectories. By convexity, the linearisation errors have tight bounds and can be treated as bounded disturbances within a robust tube MPC framework.
  We describe a systematic, data-driven method for computing DC model representations using deep learning neural networks with a special convex structure, and explain how the resulting MPC optimisation can be solved using convex programming.
  For the problem of maximising product formation in a cultivation with uncertain model parameters, we design a controller that ensures robust constraint satisfaction and allows online estimation of unknown model parameters.
  Our results indicate that this method is a promising solution for computationally tractable, robust MPC of bioprocesses.
\end{abstract}

% %%Graphical abstract
% \begin{graphicalabstract}
% %\includegraphics{grabs}
% \end{graphicalabstract}

%%Research highlights
% \begin{highlights}
% \item
%   The paper uses difference of convex (DC) function decompositions of nonlinear dynamics and successive linearisation around predicted trajectories to derive a computationally efficient convex framework for robust tube nonlinear model predictive control (MPC).

% \item
%   The paper uses deep learning neural networks to learn the model dynamics in DC form, thus providing a systematic data-driven method of computing DC model decompositions suitable for use in MPC.

% \item
%   The tube MPC algorithm is implemented using simplexes to reduce computation and improve the scalability of the approach.

% \item
%   We propose a robust adaptive tube MPC algorithm based on data-driven DC models and set membership estimation.

% \item
%   The paper describes the application of the robust adaptive MPC strategy to a bioreactor case-study.

% \end{highlights}

\begin{keyword}
Nonlinear model predictive control\sep robust adaptive control \sep data-driven control \sep convex optimisation \sep bioprocesses
%% keywords here, in the form: keyword \sep keyword
%
%% PACS codes here, in the form: \PACS code \sep code
%
%% MSC codes here, in the form: \MSC code \sep code
%% or \MSC[2008] code \sep code (2000 is the default)
%
\end{keyword}

\end{frontmatter}

%% \linenumbers

%% main text
\section{Introduction}
\label{sec:introduction}
% \subsection{Rapid bioprocess development}
The exploding demand for cost-effective production of biologic drugs and sustainable production of goods has intensified the need for rapid development and control of bioprocesses \citep{rathore21:survey}. This has motivated the development of robotic parallel experimentation facilities capable of performing sophisticated cultivations (e.g.~fed-batch, continuous) in high throughput. Some commercially available solutions are the robolector from m2p-labs, bioREACTOR 48 from 2mag, ambr250 from Sartorius, and bioXplorer from H.EL.
In the early stages of a project, development is characterised by limited process information and a large experimental design space, posing significant challenges for operating large numbers of bioreactors in parallel \citep{cruz17:parallel}. Adaptive and robust Model Predictive Control (MPC) have been shown to significantly improve process performance and experimental efficiency, even in cases of large uncertainty in
model outputs~\citep{krausch22,rolf20:ampc,tuveri23:mhe}. However, applications have been restricted to relatively stable process conditions. For example, \citet{kager20} were able to increase total product formation in a fungal process, but their approach is limited to the nominal case. The approach of~\citet{mowbray22} used neural networks (NN) to handle uncertainties but required heavy offline training.

%Increasing demand for cost-effective production of biologic drugs and sustainable biomaterials intensifies the need for rapid bioprocess development. This is particularly true in the early stages of a project, characterised by limited process information and a broad spectrum of potential optimal conditions. Advanced control approaches such as  MPC coupled with online parameter estimation have proven successful even when incomplete process information is available \citep{krausch22}, but these have been restricted to relatively stable process conditions. For example, \citet{kager20} were able to increase total product formation in a fungal process, but their approach is limited to the nominal case. The approach of~\citet{mowbray22} used neural networks (NN) to handle uncertainties but required heavy offline training. 

%\subsection{Tube-based MPC using differences of convex functions}
A popular approach in advanced control to deal with uncertain dynamic systems is tube-based MPC (TMPC). Robust nonlinear MPC requires online solution of nonconvex optimisation problems, which can be computationally expensive. A common strategy for applying TMPC to nonlinear systems is to treat the effects of nonlinearity as bounded disturbances and compute successive linear approximations around predicted trajectories~\citep{cannon11}. These approaches, nevertheless, rely on conservative estimates of linearisation error and can give poor performance~\citep{yu13}.

Recent studies have shown that tighter bounds on linearisation errors can be achieved if the problem can be expressed as a difference of convex functions~\citep{doff-sotta22,buerger24}. This is based on the observation that the necessarily convex linearisation error is maximum at the boundary of the set on which it is evaluated. Tight bounds can thus be derived and incorporated in a robust TMPC formulation. Moreover, the difference of convex functions (DC) structure of the dynamics is attractive as it results in a sequence of convex programs that can be solved with predictable computational effort.

Although any twice continuously differentiable function can be expressed in DC form, finding DC representations can be difficult. To address this issue, we use deep learning NNs with kernel weights constrained to non-negative values and convex, non-decreasing activation functions such as rectified linear units (ReLU, $\sigma(x) = \max(0,x)$), resulting in an input-convex neural network (ICNN) \citep{amos17}. The outputs of two ICNNs can be subtracted to learn dynamics in DC form~\citep{sankaranarayanan22}. Moreover, in the context of TMPC, the tube parameterisation influences the computational complexity of the optimisation problem. \citet{doff-sotta22} and \citet{buerger24}, define tube cross sections using elementwise bounds, but this requires $O(2^{n_x})$ 
inequality constraints in the MPC optimisation problem ($n_x$ denotes the number of model states),
% a number of inequality constraints in the MPC optimisation problem that depends exponentially on the number of model states,
causing significant computational burden for large systems.
Simplex tubes, for which the dependence is $O(n_x)$,
% Simplex tubes, for which the dependence is linear,
provide a computationally efficient alternative. 
This contribution describes a robust TMPC algorithm that leverages deep learning NNs to determine the process dynamics in DC form, implements simplex tubes, enables parameter estimation while providing robustness to parameter uncertainty,
and optimises product formation in a case study of a fed-batch bioreactor for the production of penicillin. 

\section{Modelling and DC approximation using deep learning}
\label{sec:modelling}

We consider a perfectly mixed isothermal fed-batch bioreactor, a popular case study example from \citep{srinivasan03}. The model states are
the cell concentration $X$ [\SI{}{\gram\per\litre}],
product concentration $P$ [\SI{}{\gram\per\litre}],
substrate concentration $S$ [\SI{}{\gram\per\litre}] and
volume $V$ [\SI{}{\litre}].
The control input is
the feed flow rate $u$ [\SI{}{\litre\per\hour}] of substrate at concentration $S_i$.
The objective of the control strategy is to maximise the product concentration $P$ after a period of $T$ hours while satisfying state and control constraints.
%constraitns on system states and control inputs.
%

The system model is given by%
\begin{equation}\label{eq:original_ct_model}
\dot{x} = \phi (x,u) + \xi
\end{equation}
where $\dot{x}=\mathrm{d}x/\mathrm{d}t$, $x = (X,S,P,V)$ is the system state and $\xi$ is an unknown bounded disturbance representing measurement noise and parameter errors not explicitly accounted for in the problem description. 
The function $\phi$ is defined for $(x,u)\in\X\times\U$, where $\X$, $\U$ are compact convex state and control constraint sets, by
% The dynamics of the system are defined for $(x,u)\in\X\times\U$ by
% \[
% \dot{x} = \phi (x,u)
% \]
% where $x = (X,S,P,V)$ is the system state, $\X, \U$ are compact state and control constraint sets, and the function $\phi$ is defined~by%
\begin{subequations}\label{eq:model}%
  \begin{align}
    \dot{X} &= \mu(S) X - \frac{u}{V_0+V} X
    \\
    \dot{S} &= -\mu(S) \frac{X}{Y_{X/S}} - v \frac{X}{Y_{P/S}} + \frac{u}{V_0+V} (S_i-S)
    \\
    \dot{P} &= vX - \frac{u}{V_0+V} P
    \\
    \dot{V} &= u
    \\
    \mu(S) &= \mu_{\max} \frac{S}{S+K_S+S^2/K_i } .
  \end{align}
\end{subequations}
Here $\mu_{\max}$ denotes the maximal growth rate (\SI{0.02}{\per\hour}),
$K_S$ is the substrate affinity constant of the cells (\SI{0.05}{\gram\per\litre}),
$K_i$ is an inhibition constant which inhibits growth at high substrate concentrations (\SI{5}{\gram\per\litre}),
$v$ is the production rate (\SI{0.004}{\litre\per\hour}),
the inlet substrate concentration is $S_i =\SI{200}{\gram\per\litre}$,
$Y_{X/S}$ is the yield coefficient of biomass per substrate ($0.4$\,$\text{g}_X$\,$\text{g}_S^{-1}$),
and $Y_{P/X}$ is the yield coefficient of product per substrate ($1.2$\,$\text{g}_P$\,$\text{g}_S^{-1}$).
The initial volume is $V_0 =\SI{120}{\litre}$.
We assume $Y_{X/S}$ and $S_i$ are unknown to the controller but lie within known intervals as a result of variability in the phenotype of the microorganism and uncertainty in feed concentration due to imperfect measurements or mixing, respectively.
% due to imperfect measurements or mixing.
%common in practice~\citep[e.g.][]{fiedler23:do-mpc}.

%
The initial conditions are
$X(0) = \SI{1}{\gram\per\litre}$,
$S(0) = \SI{0.5}{\gram\per\litre}$,
$P(0) = \SI{0}{\gram\per\litre}$,
and
$V(0) =\SI{0}{\litre}$, and the constraint sets are
\begin{equation}\label{eq:XU}
\X = \bigl\{ x : 0 \leq x \leq (3.7,\, 1,\, 3,\, 5) \bigr\}, \ \  \U  = \bigl\{ u : 0 \leq u \leq 0.2 \bigr\}.
\end{equation}

A feedforward deep learning NN framework was used to approximate the dynamics (\ref{eq:model}a)-(\ref{eq:model}e) as a difference of convex functions by subtracting the outputs of two feedforward ICNN subnetworks. An ICNN with $L$ layers is characterised by a parameter set
$\theta = \{ G_{1:L-1}, \, H_{0:L-1}, \, b_{0:L-1} \}$
and an input-output map $z_L= f(y; \theta)$, defined for $l = 0,\ldots ,L-1$ (with $(G_0, z_0) \equiv 0$) by
\begin{equation}\label{eq:NNrecursion}
z_{l+1} = \sigma(G_l z_l + H_l y + b_l) ,
\end{equation}
where $y$ is the input, $z_l$ is the layer activation, $G_l$ are positively constrained kernel weights ($\{G_l\}_{ij} \geq 0$, for all $i,j,l$),
% $l\in \{1,\ldots,L-1\}$,
$H_l$ are input passthrough weights, $b_l$ are bias terms and $\sigma(\cdot)$ is a componentwise convex, nondecreasing activation function. 
We use ReLU activation functions: $\sigma(\xi) = \max(\xi,0)$. 
% componentwise convex and nondecreasing (ReLU) activation function.
Each network layer thus consists of the composition of a linear (and hence convex) function with a nondecreasing convex function, so $z_{L}= f(y;\theta)$ is necessarily a convex map from $y$ to $z_L$. We define $y=(x,u)$ and use $z_{L}$ to approximate the map $\phi(x,u)$, where $x=(X,S,P,V)$ and $u$ are the state and control input of the system (\ref{eq:original_ct_model}). Two ICNNs whose outputs are subtracted are trained simultaneously to learn the dynamics of (\ref{eq:original_ct_model}) as a difference of convex functions:
%$f_1$, $f_2$:
\begin{equation}\label{eq:dc_ct_model}
  \dot{x} = f_1(x,u) - f_2(x,u) + \zeta,
\end{equation}
where $\zeta$ is an unknown but bounded disturbance that accounts for the approximation errors in $f_1,f_2$ and the disturbance $\xi$ in~(\ref{eq:original_ct_model}).
Each of the two ICNNs consists of an input layer, two hidden layers with $32$ nodes each, and an output layer. Thus $\phi(x,u)$ in~(\ref{eq:original_ct_model}) is approximated by $G_{2,1} z_{2,1}- G_{2,2} z_{2,2}$ where $z_{2,i}$ for $i=1$ or $2$ has the form of~(\ref{eq:NNrecursion}) with $L=2$. 

\begin{figure}[t]
\centering
\includegraphics[scale=0.44, trim=30mm 14mm 12mm 18mm, clip]{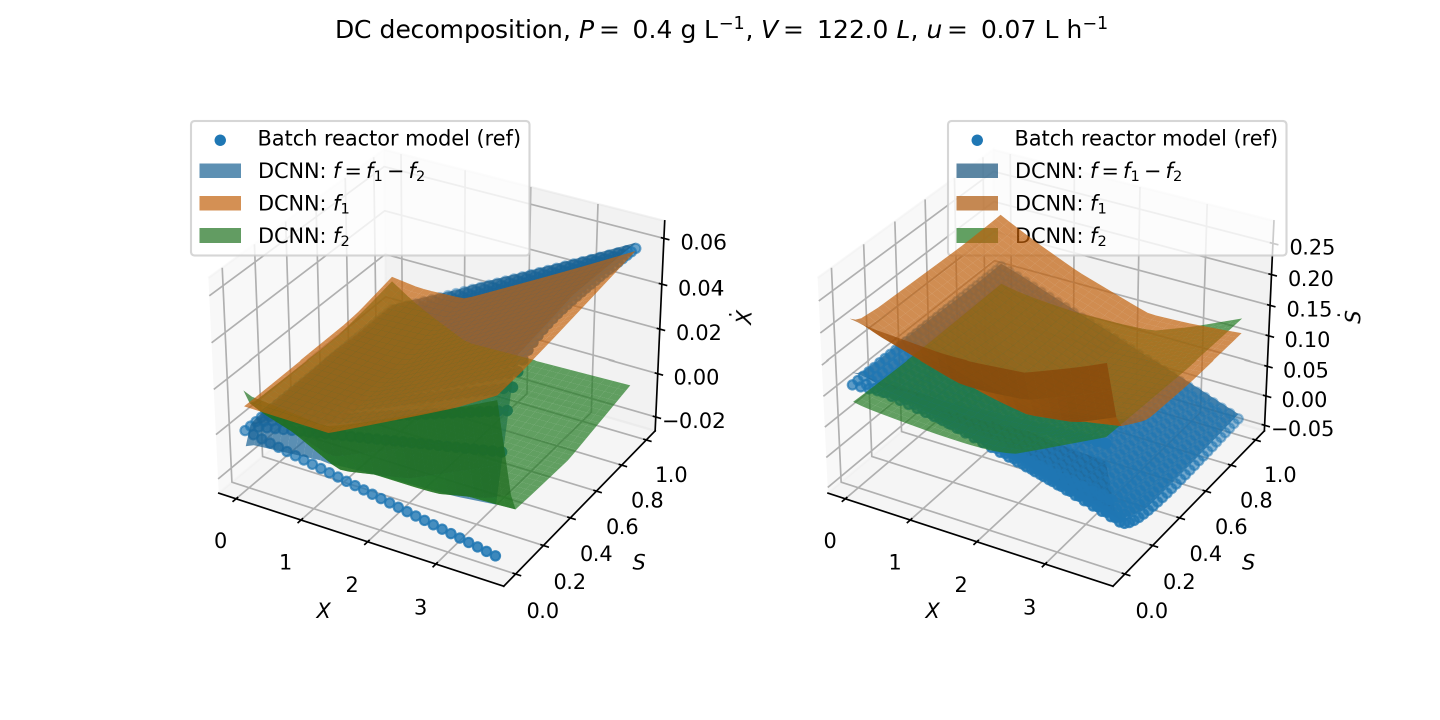}%
\vspace{-2mm}%
\caption{DC decomposition: the actual model evaluated at $30\times 30$ sample points (blue dots), the convex functions $f_1$ and $f_2$ (orange and green surfaces), and the DC decomposition $f=f_1-f_2$ (blue surfaces) approximating $\dot{X}$ (left), and $\dot{S}$ (right), as functions of the concentrations $X$ and $S$, with product concentration $P= $\,\SI{0.4}{\gram\per\litre}, volume $V=$\,\SI{122}{\litre}, and feed rate $u =$\,\SI{0.07}{\litre\per\hour}.
 \label{fig:training}}
\end{figure}

The network was implemented using Keras~\citep{chollet15:keras} and trained over $250$ epochs with the RMSProp optimiser using $10^5$ random samples of
the right-hand sides of (\ref{eq:model}a)-(\ref{eq:model}d),
%$\phi(x,u)$ in (\ref{eq:model}),
which were divided into $80\%$ training and $20\%$ validation sets.
To enable evaluation of
% $\phi(x,u)$ at
each sampling point, the unknown parameters $Y_{X/S}$, $S_i$ were set to nominal values, denoted 
$\bar{Y}_{X/S}$, $\bar{S}_i$.
Convexity of the models was evaluated by checking the Hessian matrices: $\nabla^2 f_i (x, u) \succeq 0$, $i=1,2$, at all points $(x,u)$ in the validation set. Figure~\ref{fig:training} shows the DC decompositions representing $\dot{X}$ and $\dot{S}$ for
given values of $P$, $V$, and $u$.
% a given state and control input pair $(x,u)$.
As illustrated, the NN provides obtain a good fit (MAE: $0.0012$) for the model (\ref{eq:model}) (blue dots and blue surface), and the DC form of the decomposition is apparent in Fig.~\ref{fig:training} (orange and green surfaces).   

\section{DC-TMPC framework with simplexes}\label{sec:dc-tmpc}

\citet{doff-sotta22} proposed a robust TMPC algorithm based on successive linearisation for DC systems. The so-called DC-TMPC algorithm capitalises on the idea that optimisation problems with inequality constraints involving DC functions can be approximated as convex problems by linearising only the concave parts of the inequality constraints. Furthermore, bounds on the errors introduced by linearising concave functions are trivial to compute and can be treated as disturbances in a robust TMPC scheme. We apply this approach to the system~(\ref{eq:dc_ct_model}) learned in DC form with additive disturbances, and we extend it by considering state tubes defined in terms of simplexes in order to reduce online computational load.

We can derive a discrete time model from~(\ref{eq:dc_ct_model}) using, for example, forward Euler approximation:
\begin{equation}\label{eq:dc_dt_model}
x_{k+1} = \delta \bigl( f_1(x_k, u_k) - f_2(x_k, u_k)\bigr) + w_k ,
\end{equation}
where $y_k$ denotes the value of variable $y$ at discrete time index $k\in\{0,1,2,\ldots\}$ and the sampling interval is $\delta=\SI{1}{hour}$. 
The disturbance $w_k$ lies in a bounded set $\W$ for all $k$, where $\W$ can be estimated using data on the disturbance affecting
%data from
the continuous time system (\ref{eq:original_ct_model}) and the error in its DC approximation in~(\ref{eq:dc_ct_model}).
%
%The disturbance $w_k$ belongs to a bounded set, $w_k\in\W$ for all $k$, where $\W$ can be estimated using data from the continuous time system (\ref{eq:original_ct_model}) and the error in its DC approximation in~(\ref{eq:dc_ct_model}).

%data on the continuous time system (\ref{eq:original_ct_model})-(\ref{eq:model}) and its DC approximation in~(\ref{eq:dc_ct_model}).

The tube MPC strategy is based on linearising $f_1$ and $f_2$ in~(\ref{eq:dc_dt_model}) around nominal predicted trajectories $\{x_k^\circ, u_k^\circ\}_{k=0:N}$ over a future horizon of $N$ discrete time steps. For a given control sequence $\{ u_k^\circ\}_{k=0:N}$,  we compute $\{x_k^\circ\}_{k=0:N}$ by applying an ODE solver to the model 
% by applying CasADi to
(\ref{eq:original_ct_model}) with nominal parameters (i.e.~$\xi = 0$, $Y_{X/S} = \bar{Y}_{X/S}$, $S_i = \bar{S}_i$).
We denote the state and input perturbations relative to $x_k^\circ$, $u_k^\circ$ as $s_k=x_k-x_k^\circ$ and $v_k=u_k-u_k^\circ$ respectively.
To reduce the effects of model uncertainty on predicted trajectories, $v_k$ is defined by a two degree-of-freedom control law: $v_k=K_k s_k+c_k$, where $\{K_k\}_{k=0:N}$ are feedback gains and $\{c_k\}_{k=0:N}$ is a feedforward control sequence computed at each discrete time step by solving a sequence of convex optimisation problems.

Bounds on the predicted trajectories of the model~(\ref{eq:dc_dt_model}) are obtained by exploiting the componentwise convexity of $f_i$, $i=1,2$, which implies that each component of $f_i$ is upper-bounded by the Jacobian linearisation of $f_i$ about any point in its domain. This provides the following bounds, which are by construction convex conditions on $\{s_k\}_{k=0:N}$ and $\{c_k\}_{k=0:N}$, 
\begin{subequations}\label{eq:s_dynamics}
  \begin{align}
    & s_{k+1} \geq s_k 
    + \delta \bigl( f_1(x_k^\circ + s_k, u_k^\circ + v_k) - f_1(x_k^\circ, u_k^\circ) \bigr)
    \\
    & \qquad\ \ - F_{2,k} s_k - B_{2,k} c_k + \delta \max_{w\in\W} w_k
    \nonumber \\
    & s_{k+1} \leq s_k
    + F_{1,k} s_k + B_{1,k} c_k
    \\
     & \qquad\ \ - \delta \bigl( f_2(x_k^\circ + s_k, u_k^\circ + v_k)  - f_2(x_k^\circ, u_k^\circ) \bigr) + \delta \min_{w\in\W} w_k
       \nonumber
    % \\
    % & F_{i,k} = \delta( A_{i,k} + B_{i,k} K) ,
    % \ A_{i,k} = \frac{\partial f_i}{\partial x} (x_k^\circ , u_k^\circ),
    % \ B_{i,k}= \frac{\partial f_i}{\partial u} (x_k^\circ, u_k^\circ ). 
  \end{align}
  \end{subequations}
where $F_{i,k} = \delta( A_{i,k} + B_{i,k} K) $, $A_{i,k} = \tfrac{\partial f_i}{\partial x} (x_k^\circ , u_k^\circ)$, $B_{i,k}= \tfrac{\partial f_i}{\partial u} (x_k^\circ, u_k^\circ )$ for $i=1,2$. 

To ensure robust constraint satisfaction and convergence of the successive linearisation scheme, we define an uncertainty tube consisting of a sequence of sets $\{\S_k\}_{k=0:N}$ bounding the state perturbation trajectories under all realisations of model uncertainty. The tube $\{\S_k\}_{k=0:N}$ may be defined in terms of componentwise bounds~\citep[e.g.][]{doff-sotta22}, but the number of constraints needed then grows exponentially with the state dimension $n_x$. To reduce this to a linear dependence, we instead parameterise $\{\S_k\}_{k=0:N}$ in terms of simplexes:
\begin{equation}\label{eq:S_def}
\S_k = \biggl\{ s \in \mathbb{R}^{n_x} :
\begin{bmatrix} {-I} \\ \underline{1}^\top \end{bmatrix} s \leq \begin{bmatrix} \alpha_k \\ \beta_k\end{bmatrix}\biggr\}
\end{equation}
where $I\in\mathbb{R}^{n_x\times n_x}$ is the identity matrix, $\underline{1}$ is a vector of ones, and the vector
$\alpha_k$ and scalar $\beta_k$ are optimisation variables.
Conditions on $\{\alpha_k,\beta_k\}_{k=0:N}$ to ensure that
the state perturbation $s_k$ belongs to $\S_k$ for all $k\in \{0,\ldots,N\}$ can be defined recursively as follows. Combining (\ref{eq:s_dynamics}) and (\ref{eq:S_def}), we require, for all $s\in \S_k$, $k=0,1,\ldots,N-1$,
\begin{subequations}\label{eq:S_dyn}%
\begin{align}
% \max_{s \in \mathbb{V}(\S_k)}⁡\Bigl\{
\alpha_{k+1} &\geq \max_{s\in\S_k} \Bigl( {-s}
    +\delta f_2 (x_k^\circ+s , u_k^\circ + K_k s + c_k )
    \\
    & \qquad\quad - \delta f_2 (x_k^\circ , u_k^\circ ) - F_{1,k} s - B_{1,k} c_k\Bigr)
    - \delta\min_{w\in\W} w
\nonumber \\
\beta_{k+1} &\geq \max_{s\in\S_k} \underline{1}^\top \Bigl( s + \delta f_1 (x_k^\circ + s , u_k^\circ + K_k s + c_k)
    \\
    & \qquad\qquad\quad - \delta f_1 (x_k^\circ,  u_k^\circ ) - F_{2,k} s - B_{2,k} c_k\Bigr) + \delta\max_{w\in \W} \underline{1}^\top w 
\nonumber 
\end{align}
\end{subequations}%
The conditions in (\ref{eq:S_dyn}a)-(\ref{eq:S_dyn}b) can be equivalently expressed in terms of $n_x+1$ convex inequalities using the set of vertices of $\S_k$, denoted $\{s_k^j\}_{j=0:n_x}$, 
which are given in terms of $\alpha_k$, $\beta_k$ as%
\begin{equation}\label{eq:simplex_vertices}
  s_k^j = \begin{cases}
    -\alpha_k & j = 0\\
    -\alpha_k + (\beta_k + \underline{1}^\top \alpha_k) \, \underline{e}_j & j =1,\ldots,n_x
    \end{cases}
  % \mathbb{V}(\S_k) = \{ -\alpha_k,\ -\alpha_k + (\beta_k + 1^\top\alpha_k) e_1 ,\ \ldots ,\ (\beta_k + 1^\top\alpha_k) e_n \} ,
\end{equation}
where $\underline{e}_j$ is the $j$th column of the identity matrix in $\mathbb{R}^{n_x\times n_x}$.
Exploiting again the convexity of $f_1$ and $f_2$, conditions (\ref{eq:S_dyn}a)-(\ref{eq:S_dyn}b) are therefore equivalent to the constraints, for $j=0,1,\ldots,n_x$:%
%$n_x+1$ constraints:
%\begin{subequations}\label{eq:S_dyn_vert}
\begin{align*}
  \alpha_{k+1} &\geq {-s_k^{j}}
                 + \delta f_2 (x_k^\circ+s^{j} , u_k^\circ + K_k s_k^{j} + c_k )- \delta f_2 (x_k^\circ , u_k^\circ )
  \nonumber \\
  & \quad - F_{1,k} s_k^{j} - B_{1,k} c_k - \delta w_{\min}
  \\
  \beta_{k+1} &\geq \underline{1}^\top \Bigl( s_k^{j}
                + \delta f_1 (x_k^\circ + s_k^j , u_k^\circ + K_k s_k^j + c_k)- \delta f_1 (x_k^\circ,  u_k^\circ )
    \nonumber \\
  &\quad\qquad - F_{2,k} s_k^j - B_{2,k} c_k\Bigr) + \delta w_{\max}
\end{align*}
% \begin{align}
%   \alpha_{k+1} &\geq {-s_k^{j}}+\begin{aligned}[t]
%     &\delta \bigl( f_2 (x_k^\circ+s_k^{j} , u_k^\circ + K_k s_k^{j} + c_k )- f_2 (x_k^\circ , u_k^\circ )\bigr)
%     \\
%     & - F_{1,k} s_k^{j} - B_{1,k} c_k
%     - w_{\min}
%     \end{aligned}
%   \\
%   \beta_{k+1} &\geq \begin{aligned}[t] & 1^\top \Bigl( s_k^{j} +
%                 \delta\bigl( f_1 (x_k^\circ + s_k^j , u_k^\circ + K_k s_k^j + c_k)- f_1 (x_k^\circ,  u_k^\circ )\bigr)\! \\
%                   &\quad - F_{2,k} s_k^j - B_{2,k} c_k \Bigr) + w_{\max}
%                   \end{aligned}
% \end{align}
% \end{subequations}
where $w_{\min}$, $w_{\max}$ are defined by
$\underline{e}_i^\top w_{\min} = \min_{w\in\W} \underline{e}_i^\top w$ for $i=1,\ldots,n_x$, $w_{\max}=\max_{w\in\W} \underline{1}^\top w$.
%is the maximum of $ 1^\top w$ over $w\in\W$.

\section{MPC algorithm for no parameter uncertainty}\label{sec:nominal_mpc}

This section considers the case in which the parameters of (\ref{eq:model}) are known to the controller but the model~(\ref{eq:original_ct_model}) is affected by unknown bounded disturbances. We therefore assume that the discrete time model~(\ref{eq:dc_dt_model}) has no parametric uncertainty (in particular $Y_{X/S}$ and $S_i$ are equal to their nominal values $\bar{Y}_{X/S}$, $\bar{S}_i$), but is subject to disturbances with known bounds $\W\ni w_k$ for all $k$.
At each sampling instant $t =0,\delta, 2\delta,\ldots$, given the current state $x(t)$ of~(\ref{eq:model}), we optimise the feedforward sequence $\{c_k\}_{k=0:N}$ and tube parameters $\{\alpha_k,\beta_k\}_{k=0:N}$ subject to (\ref{eq:S_dyn}a,b) and
$x_k\in \X$, $u_k \in \U$ for all $k$
% $x_k\in \X=\{x: 0 \leq x\leq x_{\max}\}$, $u_k \in \U=\{u: 0 \leq u\leq u_{\max}\}$, for all $k$.
by solving the following convex problem
\begin{subequations}\label{opt:mpc}
\begin{equation}
% \Bigl(\{c_k^\ast,\, \alpha_k^\ast,\, \beta_k^\ast\}_{k=0:N}\Bigr) = \arg\max
% \sum_{k=0}^{N+N_t} (l_{x,k} - \gamma l_{u,k})
\maximize_{c_k, \, \alpha_k, \, \beta_k,\, k = 0,\ldots,N}
\sum_{k=0}^{N} (l_{x,k} - \gamma l_{u,k})
\end{equation}
subject to, for all $k\in\{0,\ldots, N-1\}$ and $j\in\{0,\ldots,n_x\}$%
\begin{align}
  % System input and state constraints
 & \hspace{-1.8em} x_k^\circ + s_k^j \in \X,  \ \ u_k^\circ + K s_k^j + c_k \in \U \\
  % \U & \ni u_k^\circ + K s_k^j + c_k \\
  % \X & \ni x_k^\circ + s_k^j \\
%  0 &\leq  u_k^\circ + K s_k^j + c_k  \leq u_{\max}
%  \\
%  0 &\leq x_k^\circ + s_k^j \leq x_{\max}
% \\
  % Tube constraints
\alpha_{k+1} &\geq -s_k^j + \delta f_2 (x_k^\circ+s_k^j , u_k^\circ + K_k s_k^j + c_k ) - \delta f_2 (x_k^\circ , u_k^\circ ) \nonumber \\
  &\quad -F_{1,k} s_k^j - B_{1,k} c_k -\delta w_{\min}
\\
  \beta_{k+1} &\geq \underline{1}^\top \Bigl( s_k^j + \delta f_1 (x_k^\circ + s_k^j , u_k^\circ + K_k s_k^j + c_k) - \delta f_1 (x_k^\circ,  u_k^\circ )  \nonumber\\
  & \qquad \quad - F_{2,k} s_k^j - B_{2,k} c_k\Bigr) + \delta w_{\max}
  \\
  % Stage cost 
l_{x,k} &\leq \underline{e}_3^\top (x_{k}^\circ + s_{k}^j)
\\
l_{u,k+1} &\geq (u_{k+1}^\circ + K s_{k+1}^j + c_{k+1} -  u_{k}^\circ - K s_{k}^j - c_{k})^2
\end{align}
and the initial constraints, for all $j\in\{0,\ldots,n_x\}$%
\begin{align}
  % l_{x,0} &\leq \underline{e}_3^\top (x_{0}^\circ + s_{0}^j)
  % \\
  l_{u,0} &\geq \bigl(u_{0}^\circ + K s_{0}^j + c_{0} -  u(t-\delta)\bigr)^2
  \\
% initial constraints
  \alpha_0 &\geq x_0^\circ - x(t)
  \\
  \beta_0 &\geq \underline{1}^\top \bigl(x(t) - x_0^\circ\bigr)
  % 0 &= \alpha_0
  % \\ 
  % 0 &\leq \beta_0
\end{align}
and terminal constraints for all $k\in\{0, \ldots, N_{\mathrm{term}}\}$, $j\in\{0,\ldots,n_x\}$%
\begin{align}
% Terminal constraints
c_{N-k} & = 0 %\hspace{51.1mm}
  \\
x_{N-k}^\circ + s_{N-k}^j & \leq x_{\mathrm{term}} .
\end{align}
\end{subequations}
A brief explanation of the problem formulation is as follows.
\begin{itemize}
\item
The objective function (\ref{opt:mpc}a) aims to maximise the product concentration via the term $l_{x,k}$, which is a lower bound on $P = \underline{e}_3^\top x$ due to (\ref{opt:mpc}e). The term $l_{u,k}$ smooths the solution by penalising changes in the control input via (\ref{opt:mpc}f,g). The scalar constant $\gamma$ balances these two competing objectives. 
\item
The state and control constraints in (\ref{opt:mpc}b) depend on the sets $\X$, $\U$ defined in~(\ref{eq:XU}). By definition all model states and control inputs must be non-negative, while the upper bounds $u \leq 0.2$ and $x \leq (3.7,\, 1,\, 3,\, 5)$ are context-specific.
\item
The tube constraints (\ref{opt:mpc}c,d) are equivalent to (\ref{eq:S_dyn}a,b) and, with $s^j_k$ defined in~(\ref{eq:simplex_vertices}), ensure that $s_k\in\S_k$ for  $1 \leq k \leq N$.
%  $k\in\{1,\ldots, N\}$.
\item
The initial constraints (\ref{opt:mpc}i,j) ensure $s_0=x(t)-x_0^\circ\in\S_0$.
\item
The terminal constraints (\ref{opt:mpc}k,l) allow for a guarantee that problem~(\ref{opt:mpc}) remains feasible for all $t > 0$ if it is feasible at $t=0$. Here $N_{\mathrm{term}}$ and $x_{\mathrm{term}}\in\X$ are chosen so that the set $\{x_0 : x_k \leq x_{\mathrm{term}}, \, k = 0,\ldots,N_{\mathrm{term}}\}$ is positively invariant.
\end{itemize}

%From experiments it is known that the bioreactor is effectively shut down if the substrate concentration $S$ remains low enough for a sufficiently long time.

If the substrate concentration $S$ remains low for an extended period, the bioreactor effectively shuts down due to the depletion of the main carbon source.
Using data from simulations of (\ref{eq:original_ct_model})-(\ref{eq:model}), we identify this condition as $S \leq \SI{0.1}{\gram\per\litre}$ for a period of 5 hours, and hence we define $x_{\mathrm{term}} = (3.7, \, 0.1,\, 3,\, 5)$ and $N_{\mathrm{term}} = 5$ in the terminal constraints (\ref{opt:mpc}k,l).

The solution, denoted $\{c^\ast_k\}_{k=0:N}$, of (\ref{opt:mpc}) is used to update the state and input trajectories $\{x^\circ_k,u^\circ_k\}_{k=0:N}$ to be employed at the next iteration by setting $s_0 \gets x(t) - x_0^\circ$ and for $k=0,1,\ldots,N$:% 
\begin{subequations}\label{eq:iter_update}%
\begin{align}
u_k^\circ &\gets u_k^\circ+ K_k s_k  + c_k^\ast 
\\
s_{k+1} & \gets \Phi (x_k^\circ, u_k^\circ) - x_{k+1}^\circ
\\
x_{k+1}^\circ &\gets \Phi (x_k^\circ, u_k^\circ)
\end{align}
\end{subequations}
where $\Phi(x_k,u_k)$ is the solution for $x_{k+1}$ obtained from~(\ref{eq:original_ct_model})  under nominal conditions ($\xi=0$, $Y_{X/S} = \bar{Y}_{X/S}$, $S_i = \bar{S}_i$) using an ODE solver (e.g.~cvodes \citep{hindmarsh05:cvodes} implemented in CasADi \citep{andersson19:casadi}).

The optimisation~(\ref{opt:mpc}) and update~(\ref{eq:iter_update}) are repeated until the control perturbation $\sum_k \|c^\ast_k\|$
and change in optimal cost
fall below given thresholds, or until a maximum number of iterations is reached. The control input is then applied to the bioreactor as
$u(\tau) =u_0^\circ$ for $\tau\in[t,t+\delta)$.
At the next sampling instant, $t+\delta$, we set $x_0^\circ = x(t+\delta)$ and redefine $\{u^\circ_k , x^\circ_k\}_{k=0:N}$ using the optimal solution at the final iteration of time $t$ by setting $s_0 \gets 0$ and%
\begin{subequations}\label{eq:sample_update}%
\begin{alignat}{2}
  u_k^\circ &\gets u_{k+1}^\circ + K_{k+1} s_k + c_{k+1}^\ast , & \ \ k &= 0,1,\ldots,N-1
\\
s_{k+1} & \gets \Phi (x_k^\circ, u_k^\circ) - x_{k+1}^\circ , &\ \ k &= 0,1,\ldots,N
\\
x_{k+1}^\circ &\gets \Phi (x_k^\circ, u_k^\circ) , &\ \ k &= 0,1,\ldots,N
\end{alignat}
\end{subequations}
with $u_{N+1}^\circ \gets u^r$, where $u^r$ is such that $x(t) \leq x_{\mathrm{term}}$ for all $t>0$ if $x(0) \leq x_{\mathrm{term}}$ and $u(t) = u^r$. With $x_{\mathrm{term}} = (3.7,\, 0.1,\, 3,\, 5)$ we find $u_r = \SI{0.01}{\litre\per\hour}$ meets this requirement for the system~(\ref{eq:original_ct_model}).

Since~(\ref{eq:iter_update}) and (\ref{eq:sample_update}) use the nominal system model to update $\{u^\circ_k , x^\circ_k\}_{k=0:N}$, it cannot be guaranteed that the updated trajectory will result in a feasible updated problem~(\ref{opt:mpc}).
An alternative update strategy that provides a guarantee of feasibility is proposed in~\citet{lishkova25:tac}. However, we use here a simpler approach, exploiting the observation that the nominal trajectory $\{u^\circ_k , x^\circ_k\}_{k=0:N}$ computed at a feasible previous iteration necessarily provides a set of linearisation points such that problem~(\ref{opt:mpc}) is feasible at the current iteration. This is a direct consequence of the tube $\{\S_k\}_{k=0:N}$ containing the state trajectories for all possible uncertainty realisations.
Thus, if~(\ref{opt:mpc}) is infeasible after the update~(\ref{eq:sample_update}), 
we instead define $\{x^\circ_k,u^\circ_k\}_{k=0:N}$ as a time-shifted version of the previous nominal trajectory:%
\begin{subequations}\label{eq:fallback_update}%
\begin{alignat}{2}
u^\circ_k &\gets u^\circ_{k+1}, & \ \ k &= 0,1,\ldots,N-1
\\
x^\circ_k &\gets x^\circ_{k+1}, &\ \ k &= 0,1,\ldots,N
\\
u^\circ_{N} &\gets u^r & &
\\
x^\circ_{N+1} &\gets \Phi (x_k^\circ, u_N^\circ) & &
\end{alignat}
\end{subequations}
% $\Phi(x,u^r) \leq x_{\mathrm{term}}$ whenever $x\leq x_{\mathrm{term}}$.

The update strategy in (\ref{eq:iter_update})-(\ref{eq:fallback_update}) ensures that~(\ref{opt:mpc}) is recursively feasible, i.e.~a feasible solution of~(\ref{opt:mpc}) is obtained at each iteration and for all $t > 0$ if (\ref{opt:mpc}) is feasible at $t=0$.
The algorithm can be initialised at  $t=0$ by setting $x_0^\circ = x(0)$, $u^\circ_k = u^r$ and $c^\ast_k= 0$ for all $k$, and computing $\{x^\circ_k\}_{k=1:N}$ using (\ref{eq:iter_update}).

By construction, the DC-TMPC strategy satisfies the constraints $x(n\delta)\in\X$, $u(n\delta)\in\U$ for $n= 0,1,\ldots$ regardless of the number of iterations performed at each time step. Hence computation can be reduced to a single
% iteration, requiring only one
solution of the convex problem~(\ref{opt:mpc}) at each time step. On the other hand, if there is no limit on the number of iterations, then for the limiting case of no model uncertainty ($\W\to\{0\}$), the DC-TMPC iteration converges to an optimal solution for the problem of maximising product concentration for the system (\ref{eq:original_ct_model})-(\ref{eq:model}) subject to $(x,u)\in\X\times\U$~\citep[see][]{lishkova25:tac}.

%The DC-TMPC strategy is summarised in Algorithm~\ref{alg:mpc}.

% \subsection{Controller properties}
% Applied in closed loop to the system~(\ref{eq:original_ct_model}), the 

\section{Adaptive MPC algorithm for parameter uncertainty}\label{sec:ampc}

This section discusses robust adaptive DC-TMPC subject to parameter uncertainty and unknown disturbances with known bounds.
Since the uncertain parameters  $S_i$ and $Y_{X/S}^{-1}$ appear linearly in the bioreactor model~(\ref{eq:model}), we can write the system model~(\ref{eq:original_ct_model}) in a form that depends linearly on an unknown parameter vector $\theta$,
%$\theta = [1/Y_{X/S} - 1/\bar{Y}_{X/S} \ \ (S_i - \bar{S_i})/100]^\top$:
\begin{equation}\label{eq:linpar_ct_model}
  \dot{x} =\phi(x,u) + \underline{e}_2 \psi(x,u)^\top \theta + \xi .
\end{equation}
Here $\phi(x,u)$ is defined by (\ref{eq:model}a)-(\ref{eq:model}e) with $Y_{X/S}, S_i$ replaced by their nominal values $\bar{Y}_{X/S}, \bar{S}_i$, and $\psi(x,u)$ and $\theta$ are the vectors
\[
  \psi (x, u ) = \begin{bmatrix} -\mu(S) X \\ 100 u/(V+V_0) \end{bmatrix} ,
  \quad
\theta = \begin{bmatrix} Y_{X/S}^{-1} - \bar{Y}_{X/S}^{-1} \\ (S_i - \bar{S_i})/100 \end{bmatrix} .
\]
Extending the approach of Sections~\ref{sec:modelling}, \ref{sec:dc-tmpc} and \ref{sec:nominal_mpc} to account for uncertainty in $\theta$, we propose a robust DC-TMPC strategy for controlling (\ref{eq:linpar_ct_model}). 
This is combined with online estimation of $\theta$ to create a robust adaptive control algorithm.
% given bounds on $\theta$ and $\xi$.
% This is combined with~\citep{lor19:auto} for estimating $\theta$ online to construct a computationally efficient robust adaptive control algorithm.

Analogously to the DC model (\ref{eq:dc_dt_model}), we construct a DC representation $\psi(x,u) = g_1(x,u) - g_2(x,u)$, where $g_i(x,u)$ for $i=1,2$ is componentwise convex in $(x,u)$, using the difference of two ICNNs. For this example we use feedforward NNs with the same structure as the ICNNs described in Section~\ref{sec:modelling}, but with $16$ nodes in each hidden layer. 
%Using Keras and RMSprop,
The ICNNs were trained simultaneously on samples of the map $\psi(x,u)$ without knowledge of the true parameter values.
Combining this model with (\ref{eq:dc_dt_model}), which was identified using the nominal parameters $\bar{Y}_{X/S}, \bar{S}_i$, we obtain
\begin{align}\label{eq:linpar_dt_dc_model}
  x_{k+1} &= \delta \bigl( f_1(x_k,u_k) - f_2 (x_k, u_k)\bigr) \nonumber \\
  &\quad + \delta \underline{e}_2 \bigl( g_1(x_k,u_k) - g_2 (x_k, u_k)\bigr)^\top \theta + w_k . 
\end{align}
A set $\W$ bounding the disturbance $w_k$ can be constructed from bounds on the disturbance affecting the system~(\ref{eq:original_ct_model}) and the approximation errors in the DC representations of $\phi$ and $\psi$.

% Initial bounds on $\theta$ can be determined from prior bounds on the unknown parameters. For example, with nominal parameter values $\bar{Y}_{X/S} = 0.45$, $\bar{S}_i = \SI{195}{\gram\per\litre}$ and bounds $Y_{X/S}^{-1}\in [1.72, 2.72]$, $S_i \in [185,205]$, we obtain
% \[
% \theta \in \Theta_0 = \{ \theta : ({-0.5}, {-0.1}) \leq  \theta \leq (0.5, 0.1) \}.
% \]
%
%
An initial bounding set $\Theta_0$ containing $\theta$ can be determined if prior  parameter bounds are available. For example, with nominal parameter values $\bar{Y}_{X/S} = 0.45$, $\bar{S}_i = \SI{195}{\gram\per\litre}$ and bounds $Y_{X/S}\in [0.333,0.5]$, $S_i \in [180,220]$, we obtain
\[
\theta \in \Theta_0 = \{ \theta : ({-0.222}, {-0.15}) \leq  \theta \leq (0.778, 0.25) \},
\]
and the actual parameter values $Y_{X/S} = 0.4$, $S_i = \SI{200}{\gram\per\litre}$ imply $\theta^\ast = (0.278,0.05)$ is the true value of $\theta$.

We use online set membership estimation~\citep{lor19:auto,lu21:ijrnc} to update bounds on $\theta$ using online state measurements and prior disturbance bounds.
At time $t=n\delta$, for $n=0,1,\ldots$, the parameter set $\Theta_n = \{\theta : H \theta \leq h_n\}$ is defined in terms of a vector $h_n$ that is updated online and a constant matrix $H$ that is chosen offline so that $\Theta_n$ has a fixed
number of vertices. 
%, to keep the number of vertices of $\Theta$ constant, $H$ is a constant matrix chosen offline, and $h_n$ is computed online.
Thus if $H = [I \ \, {-I}]^\top$ so that $h_n$ determines upper and lower bounds on the elements of $\theta$, then $\Theta_n$ has $4$ vertices.
Let $D_n = D(x_n,u_n) = \delta \underline{e}_2 \bigl( g_1(x_n,u_n) - g_2 (x_n, u_n)\bigr)^\top$ and $d_n = d(x_n,u_n) = \delta \bigl( f_1(x_n,u_n) - f_2 (x_n, u_n)\bigr)$.
Then~(\ref{eq:linpar_dt_dc_model}) can be written as $x_{n+1} = D_n \theta + d_n + w_n$, and at time step $n+1$, $\theta$ must belong to an unfalsified parameter set given by
$\{\theta :  x_{n+1} - d_n - D_n \theta \in \W\}$.
Hence, using observations of $x_n$ over a window of $N_{\mathrm{est}}$ discrete time steps, we update $\Theta_n$ by solving a linear program:
\[
  h_{n,i} =   \begin{aligned}[t] & \max_{\theta\in\Theta_{n-1}} \, H_i \, \theta 
  \\
  & \text{subject to} \ x_{n-k+1} - d_{n-k} - D_{n-k}\theta \in \W, \ \forall k \in\{1,\ldots,N_{\mathrm{est}}\}
  \end{aligned}
\]
for each element $i$ of $h_n$, where $h_{n,i}$, $H_i$ are the $i$th element and $i$th row of $h_n$, $H$. This update rule ensures that $\Theta_n\subseteq \cdots \subseteq \Theta_0$ and $\Theta_n$ converges asymptotically to the true parameter vector $\theta^\ast$ if the control input is persistently exciting~\citep{lu21:ijrnc}.

We also compute  a point estimate $\hat{\theta}_n$ online using Least Mean Squares (LMS)~\citep{lor19:auto}.
%
% Let  $\hat{x}_{1|n} = D_n\hat{\theta}_n + d_n$ denote the single step prediction of~(\ref{eq:linpar_dt_dc_model}) with
Given the nominal 1-step ahead prediction from~(\ref{eq:linpar_dt_dc_model}),
% with $\theta=\hat{\theta}_n$
$\hat{x}_{1|n} = D_n\hat{\theta}_n + d_n$,
and the actual system state $x_{n+1}$, the estimate is updated at discrete time $n+1$ via
\[
\hat{\theta}_{n+1} = \Pi_{\Theta_{n+1}}\bigl[ \hat{\theta}_n + \hat{\mu} D_n^\top (x_{n+1} - \hat{x}_{1|n}) \bigr] 
\]
where $\Pi_\Theta[\cdot]$ denotes the (Euclidean) projection onto $\Theta$, and $\hat{\mu}$ is an update gain satisfying $\tfrac{1}{\hat{\mu}} > \| D(x,u) \|^2$ for all $(x,u)\in\X\times\U$.

The parameter set estimate $\Theta_n$ can be used to make the DC-TMPC optimisation~(\ref{opt:mpc}) robust to parameter uncertainty, while the point estimate $\hat{\theta}_n$ can improve the accuracy of the nominal predicted trajectories $\{x^\circ_k,u^\circ_k\}_{k=0:N}$ at sampling instant $t=n\delta$. These modifications result in an adaptive DC-MPC algorithm with improved robustness to model uncertainty and potentially better performance than the strategy outlined in Section~\ref{sec:nominal_mpc},
albeit at the expense of increased computation.
%although this comes at the expense of increased computation.

In particular, when computing the nominal predicted trajectories $\{x^\circ_k,u^\circ_k\}_{k=0:N}$ 
the map $\Phi(x,u)$ in (\ref{eq:iter_update})-(\ref{eq:fallback_update}) is determined by
solving~(\ref{eq:original_ct_model}) with $\xi=0$ and the current parameter estimates
%$Y_{X/S}, S_i$ determined by $\hat{\theta}_n$ 
\[
  Y_{X/S} = (\hat{\theta}_{n,1} + \bar{Y}_{X/S}^{-1})^{-1}, \quad S_i = 100\, \hat{\theta}_{n,2} + \bar{S}_i ,
\]
where $\hat{\theta}_n = (\hat{\theta}_{n,1},\hat{\theta}_{n,2})$. In addition, the tube $\{x^\circ_k + \S_k\}_{k=0:N}$ necessarily contains the system state under all model possible uncertainty realisations when constraints (\ref{opt:mpc}c,d) are replaced by%
% enforcing (\ref{eq:S_dyn}a,b)
\begin{subequations}\label{eq:robust_tube_constraints}
\begin{align}
  % Tube constraints
  \alpha_{k+1} &\geq -s_k^j + \delta f_2 (x_k^\circ+s_k^j , u_k^\circ + K_k s_k^j + c_k )  - \delta f_2 (x_k^\circ , u_k^\circ ) \nonumber \\
  &\quad - ( F_{1,k} s_k^j + B_{1,k} c_k )  - \delta w_{\min} \nonumber\\
  &\quad + \underline{e}_2 {\theta^m_+}^\top \bigl( \delta  g_2 (x_k^\circ+s_k^j , u_k^\circ + K_k s_k^j +  c_k ) - \delta g_2 (x_k^\circ , u_k^\circ )\bigr) \nonumber \\
  &\quad - \underline{e}_2 {\theta^m_+}^\top (G_{1,k} s_k^j + B^g_{1,k} c_k) \nonumber \\
  &\quad - \underline{e}_2 {\theta^m_-}^\top \bigl( \delta  g_1 (x_k^\circ+s_k^j , u_k^\circ + K_k s_k^j + c_k ) - \delta g_1 (x_k^\circ , u_k^\circ )\bigr) \nonumber \\
  &\quad + \underline{e}_2 {\theta^m_-}^\top (G_{2,k} s_k^j + B^g_{2,k} c_k)
\\
  \beta_{k+1} &\geq \underline{1}^\top \bigl( s_k^j + \delta f_1 (x_k^\circ + s_k^j , u_k^\circ + K_k s_k^j + c_k) - \delta f_1 (x_k^\circ,  u_k^\circ ) \bigr) \nonumber\\
  & \quad - \underline{1}^\top  (F_{2,k} s_k^j + B_{2,k} c_k) + \delta w_{\max} \nonumber\\
  &\quad  + {\theta^m_+}^\top \bigl( \delta  g_1 (x_k^\circ+s_k^j , u_k^\circ + K_k s_k^j + c_k ) - \delta g_1 (x_k^\circ , u_k^\circ )\bigr) \nonumber \\
  &\quad - {\theta^m_+}^\top (G_{2,k} s_k^j + B^g_{2,k} c_k) \nonumber  \\
  &\quad - {\theta^m_-}^\top \bigl( \delta  g_2 (x_k^\circ+s_k^j , u_k^\circ + K_k s_k^j +  c_k ) - \delta g_2 (x_k^\circ , u_k^\circ )\bigr) \nonumber \\
  &\quad + {\theta^m_-}^\top (G_{1,k} s_k^j + B^g_{1,k} c_k) .
\end{align}
\end{subequations}
Here $G_{i,k} = \delta( A^g_{i,k} + B^g_{i,k} K) $, $A^g_{i,k} = \tfrac{\partial g_i}{\partial x} (x_k^\circ , u_k^\circ)$, $B^g_{i,k}= \tfrac{\partial g_i}{\partial u} (x_k^\circ, u_k^\circ )$ for $i=1,2$, and $\{\theta^m\}_{m=1:4}$ are the vertices of $\Theta_n$, with
\[
\theta^m_+ = \max\{\theta^m, 0\}, \quad \theta^m_- = \min \{ \theta^m , 0\}
\]
denoting the non-negative and non-positive components of $\theta^m$.

Problem~(\ref{opt:mpc}), with (\ref{eq:robust_tube_constraints}a,b) replacing (\ref{opt:mpc}c,d), and with $\{x^\circ_k,u^\circ_k\}_{k=0:N}$ updated via  (\ref{eq:iter_update})-(\ref{eq:fallback_update}), is necessarily recursively feasible since $\Theta_{n+1}\subseteq\Theta_{n}$ for all $n\geq 0$. In addition, the closed loop adaptive DC-TMPC strategy ensures satisfaction of the constraints $\bigl(x(t),u(t)\bigr)\in\X\times\U$ at all sampling instants $t=n\delta$, $n=0,1,\ldots$. These properties hold regardless of the number of iterations performed at each time step, allowing computationally efficient implementations of the controller. Convergence of the parameter set  and pointwise estimates is not guaranteed for the DC-TMPC algorithm as stated here, although techniques for ensuring persistence of excitation~\citep[e.g.][]{lu23:auto}
could be incorporated into the approach.

\section{Simulation results and discussion}

The control algorithm was implemented in simulations of the fed-batch bioreactor described in Section~\ref{sec:modelling}.
The control objective is to maximise product concentration over a period of $T=\SI{100}{\hour}$ subject to constraints~(\ref{eq:XU}).
Two pairs of ICNNs were trained as described in Sections~\ref{sec:modelling} and \ref{sec:ampc} to provide a DC model representation, and the DC-TMPC algorithm was formulated as described in Sections~\ref{sec:nominal_mpc} and~\ref{sec:ampc}.
%
%to provide robustness to disturbances and parameter uncertainty, and to allow online parameter estimation.
%
The online optimisation (\ref{opt:mpc}), (\ref{eq:robust_tube_constraints}a,b) is a (convex) second order cone program (SOCP) since the constraint sets $\X$, $\U$ are polytopic and the ICNNs employ ReLU (piecewise linear) activation functions. We solve this problem using cvxpy~\citep{cvxpy} to interface the gurobi solver~\citep{Gurobi}. Model linearisation was performed using TensorFlow and Keras. 

We use a prediction horison of $N=25$ time steps and cost weight $\gamma = 0.1$ in (\ref{opt:mpc}a). Based on the ICNN validation data, the disturbance bounds are given by
\[
  \W=\bigl\{w : \lvert w \rvert \leq (10^{-2}, \, 10^{-3}, \, 10^{-3}, \, 10^{-2}) \bigr\}.
\]
The gains $\{K_k\}_{k=0:N}$ were found by solving the linear quadratic regulator problem defined by the time-varying linearised model around nominal predicted trajectories with $Q = 10^{-5} I$, $R = 1$.

For the DC-TMPC algorithm
of Section~\ref{sec:nominal_mpc}
without parameter uncertainty, a single iteration requires on average \SI{1.0}{\second} and at most \SI{1.2}{\second} (Apple M3 Pro, 36 GB memory). When parameter uncertainty is included as described in Section~\ref{sec:ampc} the computation time increases to \SI{12}{\second} on average and maximum \SI{15}{\second}.

\begin{figure}[ht]
\centerline{\includegraphics[scale=0.47, trim=3mm 4mm 1mm 8.5mm, clip]{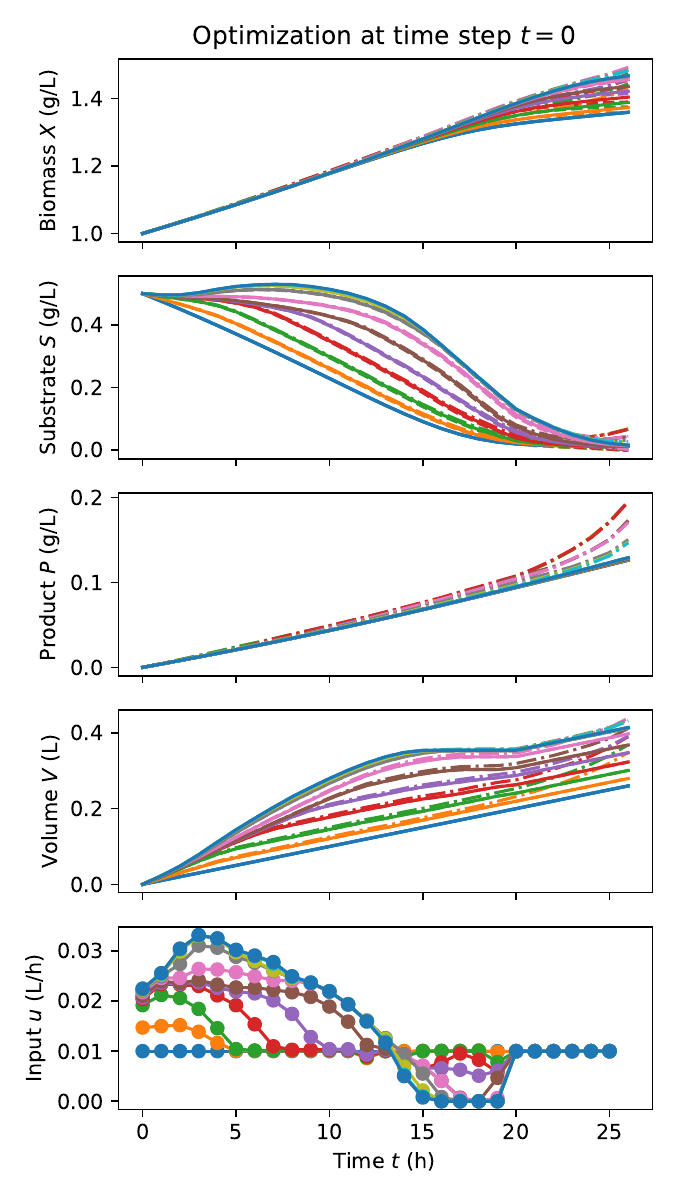}}%
\vspace{-3mm}%
\caption{Predicted state and control trajectories at initial time, $t =0$. Solid lines: nominal predicted trajectories $\{x^\circ_k,u^\circ_k\}_{k=0:N}$, dashed/dash-dotted lines: lower/upper bounds on predicted states within $\{\S_k\}_{k=0:N}$.}  
\vspace{-4mm}%
\label{fig:predicted_trajectories}
\end{figure}

Predicted trajectories computed at the first time step ($t=0$) are shown in Fig.~\ref{fig:predicted_trajectories}. Convergence was reached after $10$ iterations (at which point the changes in optimal cost and optimal solution for $\{c_k\}_{k=0:N}$ were below the chosen thresholds of $10^{-3}$ and $10^{-4}$ respectively).
% (this was the limit chosen for the maximum number of iterations)
The nominal trajectories (solid lines) can be seen to lie within the bounds on future states provided by the predicted tubes (dashed and dash-dotted lines), which is expected since the tubes provide robust bounds on the state trajectories despite the presence of uncertain parameters and unknown disturbances. The effect of the terminal constraint, requiring $S_k \leq 0.01$ for $20 \leq k \leq 25$, can also be seen.

\begin{figure}[h!]
\centerline{\includegraphics[scale=0.4, trim=1mm 4mm 1mm 2mm, clip]{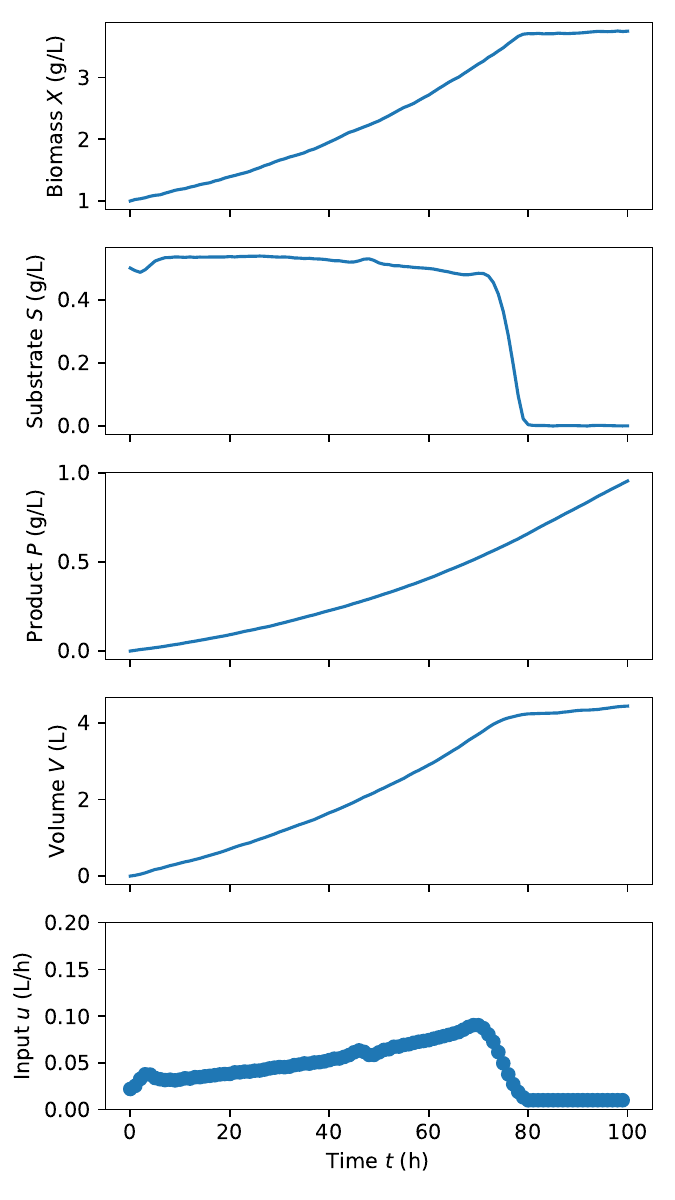}\includegraphics[scale=0.232, trim=5.9mm 22mm 8mm 20mm, clip]{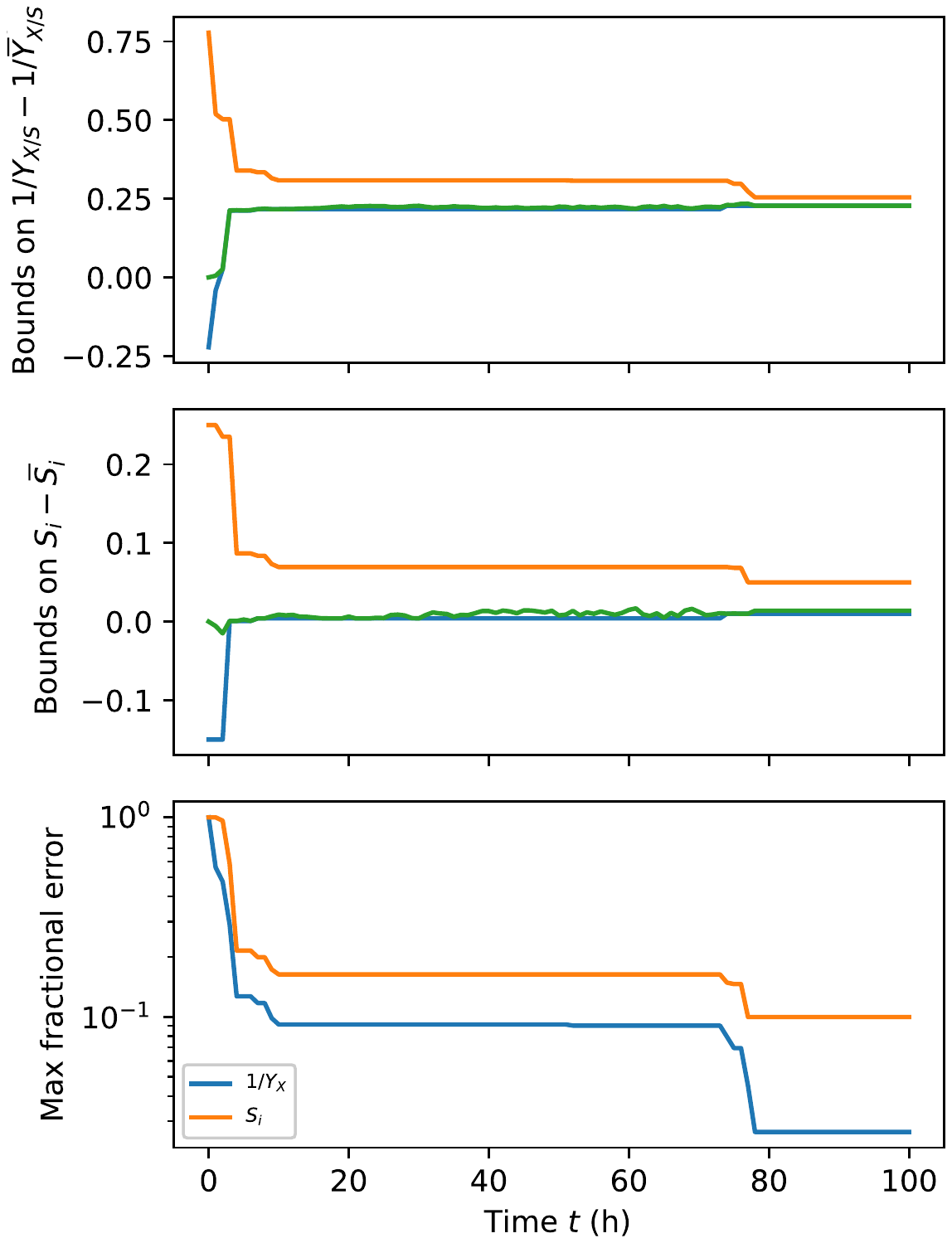}}%
\vspace{-3mm}%
\caption{Closed loop state and control trajectories (left) and parameter estimates (right).}
\label{fig:cl_trajectories}
\end{figure}

\begin{figure}[h]
\centerline{\includegraphics[scale=0.4, trim=3mm 4mm 1mm 1mm, clip]{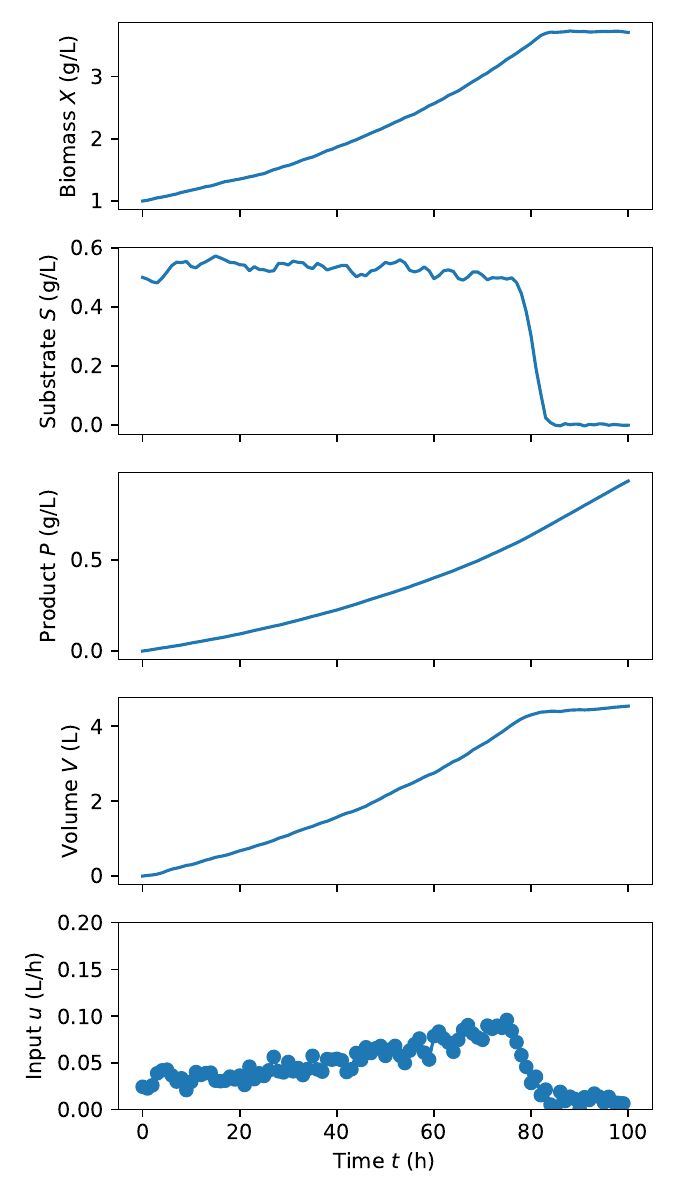}\includegraphics[scale=0.232, trim=5.9mm 22mm 8mm 20mm, clip]{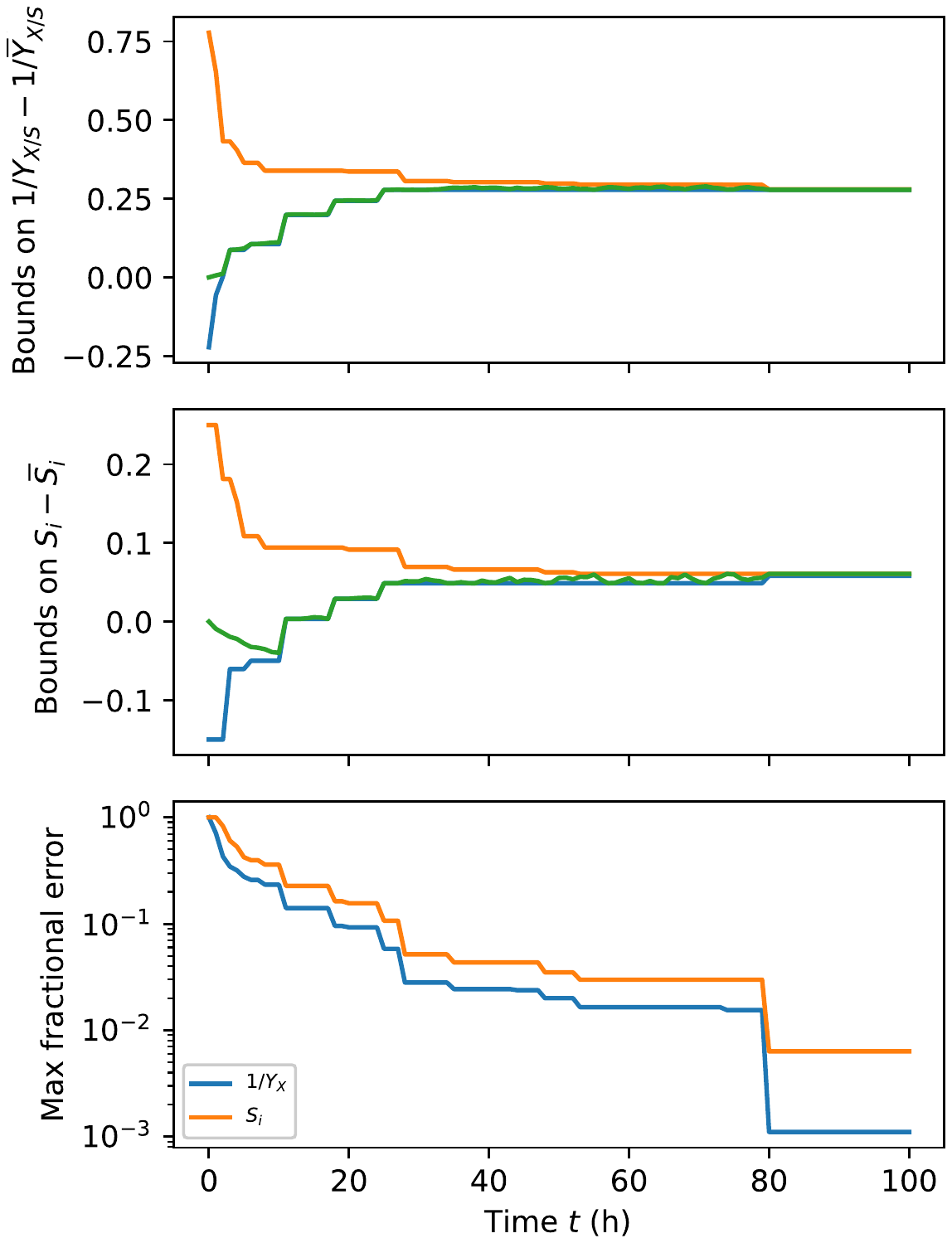}}%
\vspace{-3mm}%
\caption{Closed loop state and control trajectories with noisy inputs (left) and parameter estimates (right).}
\label{fig:cl_trajectories_dither}
\end{figure}

The closed loop state and control trajectories over $100$ time steps are shown in Fig.~\ref{fig:cl_trajectories}. The constraints on biomass ($X \leq 3.7$) and substrate concentration $S \geq 0$ are active towards the end of the simulation. Figure~\ref{fig:cl_trajectories} also
%Figure~\ref{fig:param_est}
shows the parameter estimate $\hat{\theta}_t$ (green lines) and the estimated parameter bounds corresponding to $\Theta_t$ (orange and blue lines). Although the parameter set shrinks over time, convergence is slow after the first $10$ time steps due to the unexciting optimal control trajectory in Fig.~\ref{fig:cl_trajectories}. The rate of convergence of parameter estimates improves when random noise (independent, uniformly distributed within bounds of $\pm 10^{-2}$) is added to the control sequence, as shown in Figure~\ref{fig:cl_trajectories_dither}.
% (Fig.~\ref{fig:param_est_dither}).
The resulting trajectories are more noisy %(Fig.~\ref{fig:cl_trajectories_dither}),
and the improvement in parameter estimate accuracy is accompanied by a reduction in final product concentration, which falls from $0.96$ in Fig.~\ref{fig:cl_trajectories} to $0.93$ in Fig.~\ref{fig:cl_trajectories_dither}.
The presented DC-TMPC algorithm improves on the multi-stage NMPC of~\citet{lucia13:case_study} for this case study in terms of computation and performance, and by enabling online parameter estimation.

% \begin{figure}[ht]
% \includegraphics[scale=0.232, trim=5.9mm 22mm 8mm 20mm, clip]{img/tmpc4_asymm_20-edited}%
% \vspace{-3mm}%
% \caption{Evolution of parameter set and point estimates.}
% \label{fig:param_est}
% \end{figure}

% \begin{figure}[h!]
% \centerline{\includegraphics[scale=0.4, trim=0mm 4mm 2mm 2mm, clip]{img/tmpc4_asymm_30_dither1e-2}}%
% \vspace{-3mm}%
% %\centerline{\includegraphics[scale=0.29, trim=5.9mm 22mm 8mm 20mm, clip]{img/tmpc4_asymm_dither1e-2-edited}}
% %\vspace{-6mm}%
% \caption{Evolution of parameter set and point estimates with noisy inputs.} 
% \label{fig:param_est_dither}
% \end{figure}

% \begin{figure}[h]
% \centerline{\includegraphics[scale=0.4, trim=3mm 4mm 1mm 1mm, clip]{img/tmpc1_asymm_30_dither1e-2}\includegraphics[scale=0.4, trim=0mm 4mm 2mm 2mm, clip]{img/tmpc4_asymm_30_dither1e-2}}%
% \vspace{-3mm}%
% \caption{Closed loop state and control trajectories with noisy inputs (left) and parameter estimates (right).}
% \label{fig:cl_trajectories_dither}
% \end{figure}

\section{Conclusion}
This study demonstrates the potential of deep learning robust tube MPC for optimising a bioprocess with uncertain model parameters and disturbances.
Our approach uses differences of convex (DC) functions to represent the system dynamics, and explains how to obtain these systematically using deep learning neural networks.
With this problem formulation, future performance can be optimised by solving a sequence of convex problems in which model linearisation errors appear as bounded disturbances.
Crucially, by convexity, the bounds on these errors are tight and the resulting controller is less conservative than classical TMPC based on successive linearisation.
The main contributions of this paper are the demonstration of a method using ReLU NNs for computing DC decompositions, the implementation of simplex tubes, the development of robust adaptive MPC leveraging DC models with online set membership and least mean squares estimators, and the application of these techniques within a bioreactor case-study.
Future work will incorporate more complex models and testing on a real parallel bioreactor system for bioprocess development of different host organisms and products.
%Future work will incorporate more complex models and testing on real-world bioprocesses with online optimisation of production in fast-growing E.~coli strains.

\vskip0.5\baselineskip
\emph{Acknowledgments:} We gratefully acknowledge the financial support of the German Federal Ministry of Education and Research (BMBF) (ref.~01DD20002A-KIWI Biolab) and the EPSRC (UKRI) Doctoral Prize scheme (ref.~EP/W524311/1).
\vspace{-3mm}

%% The Appendices part is started with the command \appendix;
%% appendix sections are then done as normal sections
%% \appendix

%% \section{}
%% \label{}

%% If you have bibdatabase file and want bibtex to generate the
%% bibitems, please use
%%
\bibliographystyle{elsarticle-harv} 
\bibliography{batch_bio_tmpc_nodoi}

% %% else use the following coding to input the bibitems directly in the
% %% TeX file.

% \begin{thebibliography}{00}

% %% \bibitem[Author(year)]{label}
% %% Text of bibliographic item

% \bibitem[ ()]{}

% \end{thebibliography}
\end{document}